\documentclass[10pt]{article}
\usepackage{graphicx}
\usepackage{amsmath}
\usepackage{amsfonts}
\usepackage{amssymb}
\begin{document}

\title{The cost of approximate controllability \\for semilinear heat equations \\in one space dimension }
\author{Kim Dang \textsc{Phung}\\{\small 17 rue L\'{e}onard Mafrand 92320 Ch\^{a}tillon, France}\qquad\\phung@cmla.ens-cachan.fr}
\date{}
\maketitle

\bigskip

\bigskip

Abstract: This note deals with the approximate controllability for the
semilinear heat equation in one space dimension. Our aim is to provide an
estimate of the cost of the control.

\bigskip


\bigskip

\section{Introduction and main result}

\bigskip

\bigskip

\noindent In this paper, we apply a successful combination of three key tools
which allows to get a measure of the cost of the approximate controllability
for semilinear heat equation. The first tool consists to get enough
information about the approximate control for the linear heat equation with a
potential depending on space-time variable. Then a fixed point method is
applied. The fixed point technique described here was previously used in
\cite{Z} to prove the exact controllability for semilinear wave equation in
one dimension. The last tool, usually used for control problem (see
\cite[p.589]{FCZ2} e.g.), consists to choose adequately the time of controllability.

\bigskip

\noindent Many results exist by now concerning the approximate controllability
for semilinear heat equation in a bounded domain $\Omega\subset\mathbb{R}^{n}%
$, $n\geq1$ when the control acts in a non-empty subdomain $\omega
\subset\Omega$, $\omega\neq\Omega$ (see \cite{FPZ},\cite{K} or \cite{FCZ2} and
references therein). In particular, it is proved in \cite{FCZ2} that for any
time $T>0$, if the system
\begin{equation}
\left\{
\begin{array}
[c]{c}%
\partial_{t}u-\Delta u+f\left(  u\right)  =h\cdot1_{\omega}\quad
\text{in}~\Omega\times\left(  0,T\right)  \text{ ,}\\
u=0\quad\text{on}~\partial\Omega\times\left(  0,T\right)  \text{ ,}\\
u\left(  \cdot,0\right)  =u_{o}\quad\text{in}~\Omega\text{ ,}%
\end{array}
\right.  \tag{1.1}\label{1.1}%
\end{equation}

\noindent with $f:\mathbb{R}\rightarrow\mathbb{R}$ locally
lipschitz-continuous, admits at least one globally defined and bounded
solution $u^{\ast}$, corresponding to the data $u_{o}^{\ast}\in L^{2}\left(
\Omega\right)  $ and $h^{\ast}\in L^{\infty}\left(  \omega\times\left(
0,T\right)  \right)  $, and further if the function $f$ satisfies%
\[
\left|  f^{\prime}\left(  s\right)  \right|  \leq c\left(  1+\left|  s\right|
^{p}\right)  \quad\text{a.e., with }p\leq1+4/n\text{ and }c>0\text{ ,}%
\]

\noindent and
\[
\underset{\left|  s\right|  \rightarrow\infty}{\lim}\frac{f\left(  s\right)
}{\left|  s\right|  \ln^{3/2}\left(  1+\left|  s\right|  \right)  }=0\text{ ,}%
\]

\noindent then for any $u_{o}\in L^{2}\left(  \Omega\right)  $, $u_{d}\in
L^{2}\left(  \Omega\right)  $ and $\varepsilon>0$, there exists a control
$h\in L^{\infty}\left(  \omega\times\left(  0,T\right)  \right)  $ such that
the solution of (\ref{1.1}) is globally defined in $\left[  0,T\right]  $ and
satisfies
\[
\left\|  u\left(  \cdot,T\right)  -u_{d}\right\|  _{L^{2}\left(
\Omega\right)  }\leq\varepsilon\text{ .}%
\]

\bigskip

\noindent However, in \cite{FCZ2}, no information was given about a measure of
the control with respect to $\varepsilon$. In this paper, we provide an
estimate of the control but under more restrictive hypothesis. Our result is

\bigskip

\noindent Theorem .- \textit{Let }$\Omega=\left(  0,1\right)  $\textit{ and
}$T>0$\textit{. Assume }$f\in C^{1}\left(  \mathbb{R}\right)  $\textit{ and}
\[
\underset{\left|  s\right|  \rightarrow\infty}{\lim}\frac{f\left(  s\right)
}{\left|  s\right|  \sqrt{\ln\left(  1+\left|  s\right|  \right)  }}=0\text{
,}%
\]

\noindent\textit{then, for any }$\left(  u_{o},u_{d}\right)  \in H_{0}%
^{1}\left(  \Omega\right)  \times H_{0}^{1}\left(  \Omega\right)  $\textit{
and any }$\varepsilon\in\left(  0,1\right]  $\textit{, there exist a control
}$h_{\varepsilon}\in L^{2}\left(  \omega\times\left(  0,T\right)  \right)
$\textit{ and a function }$u=u\left(  x,t\right)  \in L^{\infty}\left(
\Omega\times\left(  0,T\right)  \right)  $\textit{ such that}
\[
\left\|  h_{\varepsilon}\right\|  _{L^{2}\left(  \omega\times\left(
0,T\right)  \right)  }\leq\exp\left(  e^{C/\varepsilon}\right)  \text{ ,}%
\]%

\[
\left\|  u\left(  \cdot,T\right)  -u_{d}\right\|  _{L^{2}\left(
\Omega\right)  }\leq\varepsilon\text{ ,}%
\]

\noindent\textit{and}%
\begin{equation}
\left\{
\begin{array}
[c]{c}%
\partial_{t}u-\partial_{xx}u+f\left(  u\right)  =h_{\varepsilon}\cdot
1_{\omega}\quad\text{\textit{in}}~\Omega\times\left(  0,T\right)  \text{ ,}\\
u=0\quad\text{\textit{on}}~\partial\Omega\times\left(  0,T\right)  \text{ ,}\\
u\left(  \cdot,0\right)  =u_{o}\quad\text{\textit{in}}~\Omega\text{ .}%
\end{array}
\right.  \tag{1.2}\label{1.2}%
\end{equation}

\noindent\textit{Here, }$C$\textit{ is a positive constant independent on
}$\varepsilon$\textit{.}

\bigskip

\bigskip

Remark .- Notice that we do not assume $f\left(  0\right)  =0$. If $f\left(
0\right)  =0$ (which correspond to the case $u^{\ast}=0$), we can use the
following control strategy to provide an estimate of the control when
$u_{o}\in L^{2}\left(  \Omega\right)  $: we divide the time interval $\left(
0,T\right)  $ in two subintervals. During the first time interval $\left(
0,T/2\right]  $, we use a null control to steer the semilinear heat equation
starting from $u_{o}$ to zero (see \cite{FCZ2}). In the second time interval
$\left(  T/2,T\right)  $, we apply the above Theorem with null initial data.

\bigskip

\noindent The rest of this note is devoted to the proof of Theorem.

\bigskip

\bigskip

\section{Proof of Theorem}

\bigskip

\bigskip

We proceed in three steps.

\bigskip

Step 1 .- Preliminary on the cost of the approximate controllability for the
linear heat equation with a potential. We first recall some results from
\cite{P} concerning the cost of the approximate controllability for the heat
equation with a potential $a=a\left(  x,t\right)  \in L^{\infty}\left(
\Omega\times\left(  0,T\right)  \right)  $. We denote $\left\|  a\right\|
_{\infty}=\left\|  a\right\|  _{L^{\infty}\left(  \Omega\times\left(
0,T\right)  \right)  }$. In the sequel, $c_{1}>1$ and $c_{2}>1$ are two
constants only depending on $\Omega$ and $\omega$. Let $T^{\prime}\in\left(
0,T\right]  $ called time of controllability of the linear system. We
introduce the operator $\mathbf{C}$ given by
\[
\mathbf{C}:\vartheta\in L^{2}\left(  \omega\times\left(  0,T^{\prime}\right)
\right)  \longrightarrow w\left(  \cdot,0\right)  \in L^{2}\left(
\Omega\right)  \text{ ,}%
\]

\noindent where $w\in C\left(  \left[  0,T^{\prime}\right]  ;H_{0}^{1}\left(
\Omega\right)  \right)  \cap W^{1,2}\left(  0,T;L^{2}\left(  \Omega\right)
\right)  $ is the solution of
\[
\left\{
\begin{array}
[c]{c}%
-\partial_{t}w-\Delta w+aw=E\vartheta\cdot1_{\left|  \omega\right.  }%
\quad\text{in~}\Omega\times\left(  0,T^{\prime}\right)  \text{ ,}\\
w=0\quad\text{on~}\partial\Omega\times\left(  0,T^{\prime}\right)  \text{ ,}\\
w\left(  \cdot,T^{\prime}\right)  =0\quad\text{in~}\Omega\text{ ,}%
\end{array}
\right.
\]

\noindent with $a\in L^{\infty}\left(  \Omega\times\left(  0,T\right)
\right)  $ and $E=\exp\left(  c_{2}\left(  1+T^{\prime}\left\|  a\right\|
_{\infty}\left(  1+e^{c_{2}T^{\prime}\left\|  a\right\|  _{\infty}^{2}%
}\right)  +\left\|  a\right\|  _{\infty}^{2/3}\right)  \right)  $. We define
$\mathcal{F}=\operatorname{Im}\mathbf{C}$ the space of exact controllability
initial data with the following norm :
\begin{equation}
\left\|  w_{o}\right\|  _{\mathcal{F}}=\inf\left\{  \left\|  \vartheta
\right\|  _{L^{2}\left(  \omega\times\left(  0,T^{\prime}\right)  \right)
}\left\backslash {}\right.  \mathbf{C}\vartheta=w_{o}\right\}  \text{ .}
\tag{2.1}\label{2.1}%
\end{equation}
Denote $\mathbf{C}^{\ast}$ the adjoint of $\mathbf{C}$. It has been proved
(see \cite{P}) that the operator $\mathbf{B}=\mathbf{CC}^{\ast}$ is
non-negative, compact and self-adjoint on $L^{2}\left(  \Omega\right)  $ which
allows us to associate the Hilbert basis with eigenfunctions $\xi_{n}$ of
$\mathbf{B}$ and eigenvalues $\mu_{n}>0$ where $\mu_{n}$ is non-increasing and
tends to zero. Furthermore, let the sets $S_{n}=\left\{  m>0\left/  {}\right.
\alpha_{n+1}<\mu_{m}\leq\alpha_{n}\right\}  $ where
\begin{equation}
\alpha_{n}=e^{\mu_{1}+e}e^{-e^{n}}\text{ ,} \tag{2.2}\label{2.2}%
\end{equation}
for all $n>0$, then each function $\phi\in L^{2}\left(  \Omega\right)  $ can
be represented in the form $\phi=\sum_{n>0}\phi_{n}$ where $\phi_{n}%
=\sum_{m\in S_{n}}\left(  \phi,\xi_{m}\right)  \xi_{m}$. Finally, let $N>0$
and $z\in H_{0}^{1}\left(  \Omega\right)  $, then we can write, in
$L^{2}\left(  \Omega\right)  $ :
\[
z=\sum_{n\leq N}z_{n}+\sum_{n>N}z_{n}\text{ with }z_{n}=\sum_{m\in S_{n}%
}\left(  z,\xi_{m}\right)  \xi_{m}\text{ ,}%
\]

\noindent with the properties
\begin{equation}%
\begin{array}
[c]{ll}%
\left\|  \sum_{n\leq N}z_{n}\right\|  _{\mathcal{F}}\leq c_{3}\frac{1}%
{\sqrt{\alpha_{N+1}}}\left\|  z\right\|  _{L^{2}\left(  \Omega\right)  }\text{
,} & \\
\left\|  \sum_{n\leq N}z_{n}-z\right\|  _{L^{2}\left(  \Omega\right)  }\leq
c_{3}\frac{D}{\ln\left(  2+\frac{1}{\sqrt{\alpha_{N+1}}}\right)  }\left\|
z\right\|  _{H_{0}^{1}\left(  \Omega\right)  }\text{ ,} &
\end{array}
\tag{2.3}\label{2.3}%
\end{equation}

\noindent for some constant $c_{3}>0$ independent on $N$, $z$, $T^{\prime}$
and $a$ and where $D=c_{1}\left(  T^{\prime}e^{c_{1}T^{\prime}\left\|
a\right\|  _{\infty}^{2}}+\frac{1}{T^{\prime}}\right)  >1$ (see \cite{P}).
Here, $\sum_{n\leq N}z_{n}\in\mathcal{F}$ and precisely%
\[%
\begin{array}
[c]{ll}%
\sum_{n\leq N}z_{n} & =\sum_{n\leq N}\sum_{m\in S_{n}}\left(  z,\xi
_{m}\right)  \xi_{m}\\
& =\mathbf{C}\left(  \sum_{n\leq N}\sum_{m\in S_{n}}\left(  z,\xi_{m}\right)
\frac{1}{\mu_{m}}\mathbf{C}^{\ast}\xi_{m}\right)  \text{ .}%
\end{array}
\]

\bigskip

\noindent On another hand, let $\chi\cdot1_{\left|  \omega\right.  }$ be the
null-control function which steers to zero at time $T^{\prime}$ the solution
of the heat equation with potential $a\left(  x,T^{\prime}-t\right)  $ and
initial data $\pi_{o}\in L^{2}\left(  \Omega\right)  $. It is known (see
\cite{FCZ1}) that
\begin{equation}
\left\|  \chi\right\|  _{L^{2}\left(  \omega\times\left(  0,T^{\prime}\right)
\right)  }\leq G\left\|  \pi_{0}\right\|  _{L^{2}\left(  \Omega\right)
}\text{ ,} \tag{2.4}\label{2.4}%
\end{equation}

\noindent where $G=\exp\left(  c_{0}\left(  1+\frac{1}{T^{\prime}}+T^{\prime
}\left\|  a\right\|  _{\infty}+\left\|  a\right\|  _{\infty}^{2/3}\right)
\right)  $ for some constant $c_{0}>0$ only depending on $\Omega$ and $\omega$.

\bigskip

\noindent Therefore, for all $T^{\prime}\in\left(  0,T\right]  $, $a\in
L^{\infty}\left(  \Omega\times\left(  0,T\right)  \right)  $, $\pi_{o}\in
L^{2}\left(  \Omega\right)  $, $z\in H_{0}^{1}\left(  \Omega\right)  $, if we
choose
\[
\ell\left(  x,T^{\prime}-t\right)  =E\sum_{n\leq N}\sum_{m\in S_{n}}\left(
z,\xi_{m}\right)  \frac{1}{\mu_{m}}\mathbf{C}^{\ast}\xi_{m}%
\]

\noindent then from (\ref{2.1}), (\ref{2.2}), (\ref{2.3}) and (\ref{2.4}), the
solution $v_{1}\in C\left(  \left[  0,T^{\prime}\right]  ;H_{0}^{1}\left(
\Omega\right)  \right)  \cap W^{1,2}\left(  0,T^{\prime};L^{2}\left(
\Omega\right)  \right)  $ of
\[
\left\{
\begin{array}
[c]{c}%
\partial_{t}v_{1}-\Delta v_{1}+a\left(  x,T^{\prime}-t\right)  v_{1}=\left(
\chi+\ell\right)  \cdot1_{\left|  \omega\right.  }\quad\text{in~}\Omega
\times\left(  0,T^{\prime}\right)  \text{ ,}\\
v_{1}=0\quad\text{on~}\partial\Omega\times\left(  0,T^{\prime}\right)  \text{
,}\\
v_{1}\left(  \cdot,0\right)  =\pi_{o}\quad\text{in~}\Omega\text{ ,}%
\end{array}
\right.
\]

\noindent satisfies
\begin{equation}
\left\|  v_{1}\left(  \cdot,T^{\prime}\right)  -z\right\|  _{L^{2}\left(
\Omega\right)  }\leq c_{4}De^{-N}\left\|  z\right\|  _{H_{0}^{1}\left(
\Omega\right)  }\text{ ,} \tag{2.5}\label{2.5}%
\end{equation}

\noindent and moreover%
\begin{equation}
\left\|  \chi+\ell\right\|  _{L^{2}\left(  \omega\times\left(  0,T^{\prime
}\right)  \right)  }\leq G\left\|  \pi_{o}\right\|  _{L^{2}\left(
\Omega\right)  }+c_{4}Ee^{e^{N}}\left\|  z\right\|  _{L^{2}\left(
\Omega\right)  }\text{ ,} \tag{2.6}\label{2.6}%
\end{equation}

\noindent for any $N\geq N_{o}$ where $N_{o}>0$ and $c_{4}\geq e^{N_{o}}$.
Clearly, the approximate-control function $\ell$ depends on $N$, $z$ and $a$
coming from $E$ and the Hilbert basis $\left(  \xi_{n},\mu_{n}\right)  $.\bigskip

\bigskip

\noindent Next, let us introduce the operator $\mathbf{S}$ given by
\[
\mathbf{S}:\lambda\in\mathbb{R}\longrightarrow v_{2}\left(  \cdot,T^{\prime
}\right)  \in H_{0}^{1}\left(  \Omega\right)  \text{ ,}%
\]

\noindent where $v_{2}\in C\left(  \left[  0,T^{\prime}\right]  ;H_{0}%
^{1}\left(  \Omega\right)  \right)  \cap W^{1,2}\left(  0,T^{\prime}%
;L^{2}\left(  \Omega\right)  \right)  $ is the unique solution of
\[
\left\{
\begin{array}
[c]{c}%
\partial_{t}v_{2}-\Delta v_{2}+a\left(  x,T^{\prime}-t\right)  v_{2}%
=\lambda\quad\text{in~}\Omega\times\left(  0,T^{\prime}\right)  \text{ ,}\\
v_{2}=0\quad\text{on~}\partial\Omega\times\left(  0,T^{\prime}\right)  \text{
,}\\
v_{2}\left(  \cdot,0\right)  =0\quad\text{in~}\Omega\text{ ,}%
\end{array}
\right.
\]

\noindent One can easily check that
\begin{equation}
\left\|  \mathbf{S}\left(  \lambda\right)  \right\|  _{H_{0}^{1}\left(
\Omega\right)  }=\left\|  \nabla v_{2}\left(  \cdot,T^{\prime}\right)
\right\|  _{L^{2}\left(  \Omega\right)  }\leq\left|  \lambda\right|
\sqrt{T^{\prime}}e^{c_{5}T^{\prime}\left\|  a\right\|  _{\infty}^{2}}\text{ ,}
\tag{2.7}\label{2.7}%
\end{equation}

\noindent for some constant $c_{5}>0$ only depending on $\Omega$ and $\omega$.

\bigskip

\noindent Consequently, for all $T^{\prime}\in\left(  0,T\right]  $, $a\in
L^{\infty}\left(  \Omega\times\left(  0,T\right)  \right)  $, $\pi_{o}\in
L^{2}\left(  \Omega\right)  $, $z_{d}\in H_{0}^{1}\left(  \Omega\right)  $, if
we choose $z=z_{d}-\mathbf{S}\left(  \lambda\right)  $ the solution
$v_{3}=v_{1}+v_{2}\in C\left(  \left[  0,T^{\prime}\right]  ;H_{0}^{1}\left(
\Omega\right)  \right)  \cap W^{1,2}\left(  0,T^{\prime};L^{2}\left(
\Omega\right)  \right)  $ of
\[
\left\{
\begin{array}
[c]{c}%
\partial_{t}v_{3}-\Delta v_{3}+a\left(  x,T^{\prime}-t\right)  v_{3}%
=\lambda+\left(  \chi+\ell\right)  \cdot1_{\left|  \omega\right.  }%
\quad\text{in~}\Omega\times\left(  0,T^{\prime}\right)  \text{ ,}\\
v_{3}=0\quad\text{on~}\partial\Omega\times\left(  0,T^{\prime}\right)  \text{
,}\\
v_{3}\left(  \cdot,0\right)  =\pi_{o}\quad\text{in~}\Omega\text{ ,}%
\end{array}
\right.
\]

\noindent satisfies, taking into account (\ref{2.5}), (\ref{2.6}) and
(\ref{2.7}),
\[
\left\|  v_{3}\left(  \cdot,T^{\prime}\right)  -z_{d}\right\|  _{L^{2}\left(
\Omega\right)  }\leq c_{4}De^{-N}\left(  \left\|  z_{d}\right\|  _{H_{0}%
^{1}\left(  \Omega\right)  }+\left|  \lambda\right|  \sqrt{T^{\prime}}%
e^{c_{5}T^{\prime}\left\|  a\right\|  _{\infty}^{2}}\right)  \text{ ,}%
\]

\noindent and
\[
\left\|  \chi+\ell\right\|  _{L^{2}\left(  \omega\times\left(  0,T^{\prime
}\right)  \right)  }\leq G\left\|  \pi_{o}\right\|  _{L^{2}\left(
\Omega\right)  }+c_{4}Ee^{e^{N}}\left(  \left\|  z_{d}\right\|  _{L^{2}\left(
\Omega\right)  }+\left|  \lambda\right|  \sqrt{T^{\prime}}e^{c_{5}T^{\prime
}\left\|  a\right\|  _{\infty}^{2}}\right)  \text{ .}%
\]

\bigskip

\noindent Finally, let $q\in L^{\infty}\left(  \Omega\times\left(  0,T\right)
\right)  $. Now, we conclude with the construction of a solution $v$ of the
heat equation with a potential and a second member and with a control acting
on the interval $\left(  T-T^{\prime},T\right)  $. Precisely, we divide the
time interval $\left(  0,T\right)  $ in two subintervals. During the first
time interval $\left(  0,T-T^{\prime}\right]  $, we let the system%
\[
\left\{
\begin{array}
[c]{c}%
\partial_{t}v-\Delta v+qv=\lambda\quad\text{in~}\Omega\times\left(
0,T-T^{\prime}\right)  \text{ ,}\\
v=0\quad\text{on~}\partial\Omega\times\left(  0,T-T^{\prime}\right)  \text{
,}\\
v\left(  \cdot,0\right)  =u_{o}\quad\text{in~}\Omega\text{ ,}%
\end{array}
\right.
\]

\noindent to evolve freely without control. In the second time interval
$\left(  T-T^{\prime},T\right)  $, we choose $a\left(  \cdot,t\right)
=q\left(  \cdot,T-t\right)  $, $\pi_{o}=v\left(  \cdot,T-T^{\prime}\right)  $
and the control function such that%
\[
\left\{
\begin{array}
[c]{c}%
\partial_{t}v-\Delta v+qv=\lambda+\left[  \left(  \chi+\ell\right)  \left(
x,T^{\prime}-T+t\right)  \right]  \cdot1_{\left|  \omega\times\left(
T-T^{\prime},T\right)  \right.  }\quad\text{in~}\Omega\times\left(
0,T\right)  \text{ ,}\\
v=0\quad\text{on~}\partial\Omega\times\left(  0,T\right)  \text{ ,}\\
v\left(  \cdot,0\right)  =u_{o}\quad\text{in~}\Omega\text{ ,}%
\end{array}
\right.
\]

\noindent satisfies
\[
\left\|  v\left(  \cdot,T\right)  -z_{d}\right\|  _{L^{2}\left(
\Omega\right)  }\leq c_{4}De^{-N}\left(  \left\|  z_{d}\right\|  _{H_{0}%
^{1}\left(  \Omega\right)  }+\left|  \lambda\right|  \sqrt{T^{\prime}}%
e^{c_{5}T^{\prime}\left\|  q\right\|  _{\infty}^{2}}\right)  \text{ ,}%
\]

\noindent and moreover%
\[
\left\|  \chi+\ell\right\|  _{L^{2}\left(  \omega\times\left(  0,T^{\prime
}\right)  \right)  }\leq G\left\|  v\left(  \cdot,T-T^{\prime}\right)
\right\|  _{L^{2}\left(  \Omega\right)  }+c_{4}Ee^{e^{N}}\left(  \left\|
z_{d}\right\|  _{L^{2}\left(  \Omega\right)  }+\left|  \lambda\right|
\sqrt{T^{\prime}}e^{c_{5}T^{\prime}\left\|  q\right\|  _{\infty}^{2}}\right)
\text{ ,}%
\]

\noindent for any $N\geq N_{o}$ where $N_{o}>0$ and $c_{4}\geq e^{N_{o}}$.
Notice that one can easily check that
\[
\left\|  v\left(  \cdot,T-T^{\prime}\right)  \right\|  _{L^{2}\left(
\Omega\right)  }\leq e^{c_{6}T\left\|  q\right\|  _{\infty}^{2}}\left(
\left\|  u_{o}\right\|  _{L^{2}\left(  \Omega\right)  }+c_{6}\left|
\lambda\right|  \sqrt{T}\right)  \text{ ,}%
\]

\noindent for some constant $c_{6}>0$ only depending on $\Omega$ and $\omega$.

\bigskip

\noindent Choosing
\[
N\leq\ln\left(  c_{4}De\frac{1+\varepsilon}{\varepsilon}\left(  1+\left\|
z_{d}\right\|  _{H_{0}^{1}\left(  \Omega\right)  }+\left|  \lambda\right|
\sqrt{T^{\prime}}e^{c_{5}T^{\prime}\left\|  q\right\|  _{\infty}^{2}}\right)
\right)  <N+1
\]

\noindent then one has
\[
\left\|  v\left(  \cdot,T\right)  -z_{d}\right\|  _{L^{2}\left(
\Omega\right)  }\leq\varepsilon\text{ ,}%
\]

\noindent and moreover,
\[%
\begin{array}
[c]{ll}%
\left\|  \chi+\ell\right\|  _{L^{2}\left(  \omega\times\left(  0,T^{\prime
}\right)  \right)  } & \leq Ge^{c_{6}T\left\|  q\right\|  _{\infty}^{2}%
}\left(  \left\|  u_{o}\right\|  _{L^{2}\left(  \Omega\right)  }+c_{6}\left|
\lambda\right|  \sqrt{T}\right) \\
& \quad+c_{4}E\exp\left(  c_{4}De\frac{1+\varepsilon}{\varepsilon}\left(
1+\left\|  z_{d}\right\|  _{H_{0}^{1}\left(  \Omega\right)  }+\left|
\lambda\right|  \sqrt{T^{\prime}}e^{c_{5}T^{\prime}\left\|  q\right\|
_{\infty}^{2}}\right)  \right) \\
& \quad\quad\cdot\left(  \left\|  z_{d}\right\|  _{L^{2}\left(  \Omega\right)
}+\left|  \lambda\right|  \sqrt{T^{\prime}}e^{c_{5}T^{\prime}\left\|
q\right\|  _{\infty}^{2}}\right)  \text{ .}%
\end{array}
\]

\bigskip

\bigskip

Step 2 .- Introduction of $g$ and choice of $T^{\prime}$. We begin to fix
$\varepsilon\in\left(  0,1\right]  $ and $\left(  u_{o},u_{d}\right)  \in
H_{0}^{1}\left(  \Omega\right)  \times H_{0}^{1}\left(  \Omega\right)  $.
Next, we introduce
\[
g\left(  s\right)  =\left|
\begin{array}
[c]{ll}%
\frac{f\left(  s\right)  -f\left(  0\right)  }{s} & \text{for }s\neq0\\
f^{\prime}\left(  0\right)  & \text{at }s=0
\end{array}
\right.
\]

\noindent which satisfies, from our hypothesis on $f$, the following
assertion
\[
\forall\delta>0\quad\exists C_{\delta}>0\quad\forall s\in\mathbb{R}%
\quad\left|  g\left(  s\right)  \right|  \leq C_{\delta}+\delta\sqrt
{\ln\left(  1+\left|  s\right|  \right)  }\text{ ,}%
\]

\noindent and consequently, for any $u\in L^{\infty}\left(  \Omega
\times\left(  0,T\right)  \right)  $, $g\left(  u\right)  \in L^{\infty
}\left(  \Omega\times\left(  0,T\right)  \right)  $ and one has%
\[
\forall\delta>0\quad\exists C_{\delta}>0\quad\left\|  g\left(  u\right)
\right\|  _{\infty}\leq C_{\delta}+\delta\sqrt{\ln\left(  1+\left\|
u\right\|  _{\infty}\right)  }\text{ .}%
\]

\noindent Hence, we easily deduce that
\begin{equation}
\forall\delta>0\quad\exists C_{\delta}>0\quad\exp\left(  \frac{1}{\delta
}\left\|  g\left(  u\right)  \right\|  _{\infty}^{2}\right)  \leq C_{\delta
}+\left\|  u\right\|  _{\infty}\text{ .} \tag{2.8}\label{2.8}%
\end{equation}

\noindent Now, we take $T^{\prime}\in\left(  0,T\right]  $ depending on
$\varepsilon$ and $\left\|  g\left(  u\right)  \right\|  _{\infty}$ as follows%
\begin{equation}
T^{\prime}=\left|
\begin{array}
[c]{ll}%
T & \text{if }\varepsilon\left\|  g\left(  u\right)  \right\|  _{\infty}%
^{2}\leq1\\
\frac{T}{\varepsilon\left\|  g\left(  u\right)  \right\|  _{\infty}^{2}} &
\text{if }\varepsilon\left\|  g\left(  u\right)  \right\|  _{\infty}^{2}>1
\end{array}
\right.  \tag{2.9}\label{2.9}%
\end{equation}

\bigskip

\bigskip

Step 3 .- The fixed point method thanks to the homotopy invariance of the
Leray-Schauder degree. In order to prove Theorem, we will apply the
homotopical version of the Leray-Schauder fixed point theorem.

\bigskip

Theorem (Leray-Schauder) .- Let $\mathcal{E}$ be a Banach space and
$\mathbf{H}:\mathcal{E}\times\left[  0,1\right]  \rightarrow\mathcal{E}$ be a
compact continuous mapping such that $\mathbf{H}\left(  u,0\right)  =0$ for
every $u\in\mathcal{E}$. If there exists a constant $K$ such that $\left\|
u\right\|  _{\mathcal{E}}<K$ for every pair $\left(  u,\sigma\right)
\in\mathcal{E}\times\left[  0,1\right]  $ satisfying $u=\mathbf{H}\left(
u,\sigma\right)  $, then the mapping $\mathbf{H}\left(  \cdot,1\right)
:\mathcal{E}\rightarrow\mathcal{E}$ has a fixed point.

\bigskip

\bigskip

\noindent We introduce the following mapping $\mathbf{H}$
\[
\mathbf{H}:\left(  u,\sigma\right)  \in L^{\infty}\left(  \Omega\times\left(
0,T\right)  \right)  \times\left[  0,1\right]  \longrightarrow\sigma y\in
L^{\infty}\left(  \Omega\times\left(  0,T\right)  \right)  \text{ ,}%
\]

\noindent where $y\in C\left(  \left[  0,T\right]  ;H_{0}^{1}\left(
\Omega\right)  \right)  \cap W^{1,2}\left(  0,T;L^{2}\left(  \Omega\right)
\right)  $ is the solution of
\[
\left\{
\begin{array}
[c]{c}%
\partial_{t}y-\Delta y+\sigma g\left(  u\right)  y=-\sigma f\left(  0\right)
+h\cdot1_{\left|  \omega\right.  }\quad\text{in~}\Omega\times\left(
0,T\right)  \text{ ,}\\
y=0\quad\text{on~}\partial\Omega\times\left(  0,T\right)  \text{ ,}\\
y\left(  \cdot,0\right)  =u_{o}\quad\text{in~}\Omega\text{ ,}%
\end{array}
\right.
\]

\noindent when the control function $h$ depends on $\left(  u,\sigma\right)  $
as follows: from $q=\sigma g\left(  u\right)  \in L^{\infty}\left(
\Omega\times\left(  0,T\right)  \right)  $, we take $a\left(  \cdot,t\right)
=q\left(  \cdot,T-t\right)  $ and generate the eigencouple $\left(  \xi
_{n},\mu_{n}\right)  $, next we choose the control function
\[
h\left(  x,T-t\right)  =\left|
\begin{array}
[c]{ll}%
0 & \text{ for }T^{\prime}\leq t<T\\
\chi\left(  x,T^{\prime}-t\right)  +E\sum_{n\leq N}\sum_{m\in S_{n}}\left(
u_{d}-\mathbf{S}\left(  -\sigma f\left(  0\right)  \right)  ,\xi_{m}\right)
\frac{1}{\mu_{m}}\mathbf{C}^{\ast}\xi_{m} & \text{ for }0<t<T^{\prime}%
\end{array}
\right.
\]

\noindent where $N\geq N_{o}$ is such that $N\leq\ln\left(  c_{4}%
De\frac{1+\varepsilon}{\varepsilon}\left(  1+\left\|  u_{d}\right\|
_{H_{0}^{1}\left(  \Omega\right)  }+\sigma\left|  f\left(  0\right)  \right|
\sqrt{T^{\prime}}e^{c_{5}T^{\prime}\left\|  \sigma g\left(  u\right)
\right\|  _{\infty}^{2}}\right)  \right)  <N+1$, then the unique solution $y$
satisfies
\[
\left\|  y\left(  \cdot,T\right)  -u_{d}\right\|  _{L^{2}\left(
\Omega\right)  }\leq\varepsilon\text{ ,}%
\]

\noindent and moreover, one has%
\begin{equation}%
\begin{array}
[c]{ll}%
\left\|  h\right\|  _{L^{2}\left(  \omega\times\left(  0,T\right)  \right)  }%
& \leq Ge^{c_{6}T\left\|  \sigma g\left(  u\right)  \right\|  _{\infty}^{2}%
}\left(  \left\|  u_{o}\right\|  _{L^{2}\left(  \Omega\right)  }+c_{6}%
\sigma\left|  f\left(  0\right)  \right|  \sqrt{T}\right) \\
& \quad+c_{4}E\exp\left(  c_{4}De\frac{1+\varepsilon}{\varepsilon}\left(
1+\left\|  u_{d}\right\|  _{H_{0}^{1}\left(  \Omega\right)  }+\sigma\left|
f\left(  0\right)  \right|  \sqrt{T^{\prime}}e^{c_{5}T^{\prime}\left\|  \sigma
g\left(  u\right)  \right\|  _{\infty}^{2}}\right)  \right) \\
& \quad\quad\cdot\left(  \left\|  u_{d}\right\|  _{L^{2}\left(  \Omega\right)
}+\sigma\left|  f\left(  0\right)  \right|  \sqrt{T^{\prime}}e^{c_{5}%
T^{\prime}\left\|  \sigma g\left(  u\right)  \right\|  _{\infty}^{2}}\right)
\text{ ,}%
\end{array}
\tag{2.10}\label{2.10}%
\end{equation}

\noindent with%
\begin{equation}
\left\{
\begin{array}
[c]{ll}%
G=\exp\left(  c_{0}\left(  1+\frac{1}{T^{\prime}}+T^{\prime}\left\|  \sigma
g\left(  u\right)  \right\|  _{\infty}+\left\|  \sigma g\left(  u\right)
\right\|  _{\infty}^{2/3}\right)  \right)  \text{ ,} & \\
D=c_{1}\left(  T^{\prime}e^{c_{1}T^{\prime}\left\|  \sigma g\left(  u\right)
\right\|  _{\infty}^{2}}+\frac{1}{T^{\prime}}\right)  >1\text{ ,} & \\
E=\exp\left(  c_{2}\left(  1+T^{\prime}\left\|  \sigma g\left(  u\right)
\right\|  _{\infty}e^{c_{2}T^{\prime}\left\|  \sigma g\left(  u\right)
\right\|  _{\infty}^{2}}+\left\|  \sigma g\left(  u\right)  \right\|
_{\infty}^{2/3}\right)  \right)  \text{ .} &
\end{array}
\right.  \tag{2.11}\label{2.11}%
\end{equation}

\noindent Clearly, the control function $h$ depends on $\varepsilon$, $u_{o}$,
$u_{d}$ and $\left(  u,\sigma\right)  $ coming from $E$ and the eigencouple
$\left(  \xi_{n},\mu_{n}\right)  $.

\bigskip

From now, we use the letter $c$ to denote a positive constant only depending
on $\Omega$ and $\omega$, whose value can change from line to line. From
(\ref{2.10}) and (\ref{2.11}), the control function is bounded as follows: for
any $\varepsilon\in\left(  0,1\right]  $, $\left(  u_{o},u_{d}\right)  \in
L^{2}\left(  \Omega\right)  \times H_{0}^{1}\left(  \Omega\right)  $ and
$T^{\prime}\in\left(  0,T\right]  $, $T>0$,
\begin{equation}%
\begin{array}
[c]{ll}%
\left\|  h\right\|  _{L^{2}\left(  \omega\times\left(  0,T\right)  \right)  }%
& \\
\leq\left(  \left\|  u_{o}\right\|  _{L^{2}\left(  \Omega\right)  }%
+\sigma\left|  f\left(  0\right)  \right|  \sqrt{T}\right)  \exp\left(
c\left(  1+T\left\|  \sigma g\left(  u\right)  \right\|  _{\infty}^{2}%
+\frac{1}{T^{\prime}}+T^{\prime}\left\|  \sigma g\left(  u\right)  \right\|
_{\infty}+\left\|  \sigma g\left(  u\right)  \right\|  _{\infty}^{2/3}\right)
\right)  & \\
\quad+\left(  \left\|  u_{d}\right\|  _{L^{2}\left(  \Omega\right)  }%
+\sigma\left|  f\left(  0\right)  \right|  \sqrt{T}\right)  \exp\left(
c\left(  1+T^{\prime}\left\|  \sigma g\left(  u\right)  \right\|  _{\infty
}^{2}+T^{\prime}\left\|  \sigma g\left(  u\right)  \right\|  _{\infty
}e^{cT^{\prime}\left\|  \sigma g\left(  u\right)  \right\|  _{\infty}^{2}%
}+\left\|  \sigma g\left(  u\right)  \right\|  _{\infty}^{2/3}\right)  \right)
& \\
\quad\quad\cdot\exp\left(  \frac{c}{\varepsilon}\left(  T^{\prime
}e^{cT^{\prime}\left\|  \sigma g\left(  u\right)  \right\|  _{\infty}^{2}%
}+\frac{1}{T^{\prime}}\right)  \left(  1+\left\|  u_{d}\right\|  _{H_{0}%
^{1}\left(  \Omega\right)  }+\sigma\left|  f\left(  0\right)  \right|
\sqrt{T^{\prime}}e^{cT^{\prime}\left\|  \sigma g\left(  u\right)  \right\|
_{\infty}^{2}}\right)  \right)  &
\end{array}
\tag{2.12}\label{2.12}%
\end{equation}

\noindent and therefore%
\begin{equation}%
\begin{array}
[c]{ll}%
\left\|  h\right\|  _{L^{2}\left(  \omega\times\left(  0,T\right)  \right)  }%
& \\
\leq\left(  \left\|  u_{o}\right\|  _{L^{2}\left(  \Omega\right)  }%
+\sigma\left|  f\left(  0\right)  \right|  \sqrt{T}\right)  \exp\left(
c\left(  1+T\left\|  \sigma g\left(  u\right)  \right\|  _{\infty}^{2}%
+\frac{1}{T^{\prime}}+T^{\prime}\left\|  \sigma g\left(  u\right)  \right\|
_{\infty}+\left\|  \sigma g\left(  u\right)  \right\|  _{\infty}^{2/3}\right)
\right)  & \\
\quad+\left(  \left\|  u_{d}\right\|  _{L^{2}\left(  \Omega\right)  }%
+\sigma\left|  f\left(  0\right)  \right|  \sqrt{T}\right)  \exp\left(
c\left(  1+T^{\prime}\left\|  \sigma g\left(  u\right)  \right\|  _{\infty
}^{2}+\sqrt{T^{\prime}}e^{cT^{\prime}\left\|  \sigma g\left(  u\right)
\right\|  _{\infty}^{2}}+\left\|  \sigma g\left(  u\right)  \right\|
_{\infty}^{2/3}\right)  \right)  & \\
\quad\quad\cdot\exp\left(  \frac{c}{\varepsilon}\left(  1+T^{\prime
2}+T^{\prime}\left(  \left\|  u_{d}\right\|  _{H_{0}^{1}\left(  \Omega\right)
}+\left|  \sigma f\left(  0\right)  \right|  ^{2}\right)  \right)
e^{cT^{\prime}\left\|  \sigma g\left(  u\right)  \right\|  _{\infty}^{2}%
}\right)  & \\
\quad\quad\cdot\exp\left(  \frac{c}{\varepsilon T^{\prime}}\left(  1+\left\|
u_{d}\right\|  _{H_{0}^{1}\left(  \Omega\right)  }+\left|  \sigma f\left(
0\right)  \right|  ^{2}\right)  \right)  \text{ .} &
\end{array}
\tag{2.13}\label{2.13}%
\end{equation}

\bigskip

\noindent The continuity and compactness property of $\mathbf{H}$ comes from
the following embedding%
\[
W^{1,2}\left(  0,T;L^{2}\left(  \Omega\right)  \right)  \cap L^{\infty}\left(
0,T;H_{0}^{1}\left(  \Omega\right)  \right)  \subset L^{\infty}\left(
\Omega\times\left(  0,T\right)  \right)
\]

\noindent which is compact in one dimension of space. It remains to prove
that
\[
\left\|  u\right\|  _{\infty}<K\text{ ,}%
\]
for every pair $\left(  u,\sigma\right)  \in L^{\infty}\left(  \Omega
\times\left(  0,T\right)  \right)  \times\left[  0,1\right]  $ satisfying
$u=\mathbf{H}\left(  u,\sigma\right)  $.

\bigskip

\noindent The solution $u$ of the nonlinear system $\mathbf{H}\left(
u,\sigma\right)  =u$ is also solution of the linear system%
\[
\left\{
\begin{array}
[c]{c}%
\partial_{t}\psi-\partial_{xx}\psi+q\left(  x,t\right)  \psi=b\left(
x,t\right)  \quad\text{in~}\Omega\times\left(  0,T\right)  \text{ ,}\\
\psi=0\quad\text{on~}\partial\Omega\times\left(  0,T\right)  \text{ ,}\\
\psi\left(  \cdot,0\right)  =\sigma u_{o}\quad\text{in~}\Omega\text{ ,}%
\end{array}
\right.
\]

\noindent by substituting $q=\sigma^{2}g\left(  u\right)  $ and $b=\sigma
\left(  -\sigma f\left(  0\right)  +h\cdot1_{\left|  \omega\right.  }\right)
$. But such solution $\psi$ satisfies, in one space dimension, the following
inequality%
\[
\left\|  \psi\right\|  _{\infty}^{2}\leq ce^{cT\left\|  q\right\|  _{\infty
}^{2}}\left(  \left\|  \sigma u_{o}\right\|  _{H_{0}^{1}\left(  \Omega\right)
}^{2}+\left\|  b\right\|  _{L^{2}\left(  \Omega\times\left(  0,T\right)
\right)  }^{2}\right)  \text{ .}%
\]

\noindent Consequently, the later inequality and (\ref{2.13}) imply that
\[%
\begin{array}
[c]{ll}%
\left\|  u\right\|  _{\infty}^{2} & \leq ce^{cT\left\|  g\left(  u\right)
\right\|  _{\infty}^{2}}\left(  \left\|  u_{o}\right\|  _{H_{0}^{1}\left(
\Omega\right)  }^{2}+\left|  f\left(  0\right)  T\right|  ^{2}+\left\|
h\right\|  _{L^{2}\left(  \Omega\times\left(  0,T\right)  \right)  }%
^{2}\right) \\
& \leq\left(  \left\|  u_{o}\right\|  _{H_{0}^{1}\left(  \Omega\right)  }%
^{2}+\left\|  u_{d}\right\|  _{L^{2}\left(  \Omega\right)  }^{2}+\left|
f\left(  0\right)  T\right|  ^{2}\right) \\
& \quad\cdot\exp\left(  c\left(  1+T\left\|  g\left(  u\right)  \right\|
_{\infty}^{2}+\frac{1}{T^{\prime}}+\sqrt{T^{\prime}}e^{cT^{\prime}\left\|
g\left(  u\right)  \right\|  _{\infty}^{2}}+T^{\prime}\left\|  g\left(
u\right)  \right\|  _{\infty}+\left\|  g\left(  u\right)  \right\|  _{\infty
}^{2/3}\right)  \right) \\
& \quad\cdot\exp\left(  \frac{c}{\varepsilon}\left(  1+T^{\prime2}+T^{\prime
}\left(  \left\|  u_{d}\right\|  _{H_{0}^{1}\left(  \Omega\right)  }+\left|
f\left(  0\right)  \right|  ^{2}\right)  \right)  e^{cT^{\prime}\left\|
g\left(  u\right)  \right\|  _{\infty}^{2}}\right) \\
& \quad\cdot\exp\left(  \frac{c}{\varepsilon T^{\prime}}\left(  1+\left\|
u_{d}\right\|  _{H_{0}^{1}\left(  \Omega\right)  }+\left|  f\left(  0\right)
\right|  ^{2}\right)  \right)  \text{ .}%
\end{array}
\]

\bigskip

\noindent First, if $\varepsilon\left\|  g\left(  u\right)  \right\|
_{\infty}^{2}\leq1$, then we easily get an uniform bound for $u$ in
$L^{\infty}\left(  \Omega\times\left(  0,T\right)  \right)  $,%
\[%
\begin{array}
[c]{ll}%
\left\|  u\right\|  _{\infty}^{2} & \leq\left(  \left\|  u_{o}\right\|
_{H_{0}^{1}\left(  \Omega\right)  }^{2}+\left\|  u_{d}\right\|  _{L^{2}\left(
\Omega\right)  }^{2}+\left|  f\left(  0\right)  T\right|  ^{2}\right) \\
& \quad\cdot\exp\left(  c\left(  1+\frac{T}{\varepsilon}+\frac{1}{T}+\sqrt
{T}e^{cT/\varepsilon}+\frac{T}{\sqrt{\varepsilon}}+\frac{1}{\varepsilon^{1/3}%
}\right)  \right) \\
& \quad\cdot\exp\left(  \frac{c}{\varepsilon}\left(  1+T^{2}+T\left(  \left\|
u_{d}\right\|  _{H_{0}^{1}\left(  \Omega\right)  }+\left|  f\left(  0\right)
\right|  ^{2}\right)  \right)  e^{cT/\varepsilon}\right) \\
& \quad\cdot\exp\left(  \frac{c}{\varepsilon T}\left(  1+\left\|
u_{d}\right\|  _{H_{0}^{1}\left(  \Omega\right)  }+\left|  f\left(  0\right)
\right|  ^{2}\right)  \right)
\end{array}
\]

\noindent which gives%
\[
\left\|  u\right\|  _{\infty}^{2}\leq\left(  \left\|  u_{o}\right\|
_{H_{0}^{1}\left(  \Omega\right)  }^{2}+\left\|  u_{d}\right\|  _{L^{2}\left(
\Omega\right)  }^{2}+\left|  f\left(  0\right)  T\right|  ^{2}\right)
\exp\left(  C_{T}\left(  1+\left\|  u_{d}\right\|  _{H_{0}^{1}\left(
\Omega\right)  }+\left|  f\left(  0\right)  \right|  ^{2}\right)
e^{cT/\varepsilon}\right)  \text{ ,}%
\]

\noindent where $C_{T}>0$ is a constant only dependent on $T$, $\Omega$ and
$\omega$.

\bigskip

\noindent Now if $\varepsilon\left\|  g\left(  u\right)  \right\|  _{\infty
}^{2}>1$ then by the choice of $T^{\prime}$ given by (\ref{2.9}), we have%
\[%
\begin{array}
[c]{ll}%
\left\|  u\right\|  _{\infty}^{2} & \leq\left(  \left\|  u_{o}\right\|
_{H_{0}^{1}\left(  \Omega\right)  }^{2}+\left\|  u_{d}\right\|  _{L^{2}\left(
\Omega\right)  }^{2}+\left|  f\left(  0\right)  T\right|  ^{2}\right) \\
& \quad\cdot\exp\left(  c\left(  1+T\left\|  g\left(  u\right)  \right\|
_{\infty}^{2}+\frac{\varepsilon\left\|  g\left(  u\right)  \right\|  _{\infty
}^{2}}{T}+\sqrt{T}e^{cT/\varepsilon}+T\left\|  g\left(  u\right)  \right\|
_{\infty}+\left\|  g\left(  u\right)  \right\|  _{\infty}^{2/3}\right)
\right) \\
& \quad\cdot\exp\left(  \frac{c}{\varepsilon}\left(  1+T^{2}+T\left(  \left\|
u_{d}\right\|  _{H_{0}^{1}\left(  \Omega\right)  }+\left|  f\left(  0\right)
\right|  ^{2}\right)  \right)  e^{cT/\varepsilon}\right) \\
& \quad\cdot\exp\left(  c\frac{\left\|  g\left(  u\right)  \right\|  _{\infty
}^{2}}{T}\left(  1+\left\|  u_{d}\right\|  _{H_{0}^{1}\left(  \Omega\right)
}+\left|  f\left(  0\right)  \right|  ^{2}\right)  \right)
\end{array}
\]

\noindent which gives%
\[%
\begin{array}
[c]{ll}%
\left\|  u\right\|  _{\infty}^{2} & \leq\left(  \left\|  u_{o}\right\|
_{H_{0}^{1}\left(  \Omega\right)  }^{2}+\left\|  u_{d}\right\|  _{L^{2}\left(
\Omega\right)  }^{2}+\left|  f\left(  0\right)  T\right|  ^{2}\right) \\
& \quad\cdot\exp\left(  c\left(  1+T+\frac{1}{T}\left(  1+\left\|
u_{d}\right\|  _{H_{0}^{1}\left(  \Omega\right)  }+\left|  f\left(  0\right)
\right|  ^{2}\right)  \right)  \left\|  g\left(  u\right)  \right\|  _{\infty
}^{2}\right) \\
& \quad\cdot\exp\left(  \frac{c}{\varepsilon}\left(  1+T^{2}+T\left(  \left\|
u_{d}\right\|  _{H_{0}^{1}\left(  \Omega\right)  }+\left|  f\left(  0\right)
\right|  ^{2}\right)  \right)  e^{cT/\varepsilon}\right)
\end{array}
\]

\noindent and finally, using (\ref{2.8}), there exists a constant $C^{\prime
}>0$ only depending on $\left(  \left\|  u_{d}\right\|  _{H_{0}^{1}\left(
\Omega\right)  }+\left|  f\left(  0\right)  \right|  ^{2}\right)  $, $T$,
$\Omega$ and $\omega$ such that
\[
\left\|  u\right\|  _{\infty}^{2}\leq\left(  \left\|  u_{o}\right\|
_{H_{0}^{1}\left(  \Omega\right)  }^{2}+\left\|  u_{d}\right\|  _{L^{2}\left(
\Omega\right)  }^{2}+\left|  f\left(  0\right)  T\right|  ^{2}\right)
\exp\left(  C_{T}\left(  1+\left\|  u_{d}\right\|  _{H_{0}^{1}\left(
\Omega\right)  }+\left|  f\left(  0\right)  \right|  ^{2}\right)
e^{cT/\varepsilon}\right)  \left(  C^{\prime}+\left\|  u\right\|  _{\infty
}\right)  \text{ ,}%
\]

\noindent where $C_{T}>0$ is a constant only dependent on $T$, $\Omega$ and
$\omega$.

\bigskip

\noindent We conclude that any solution $\left(  u,\sigma\right)  \in
L^{\infty}\left(  \Omega\times\left(  0,T\right)  \right)  \times\left[
0,1\right]  $ of $u=\mathbf{H}\left(  u,\sigma\right)  $ satisfies the
following estimate: there is a constant $C>0$ independent of $\left(
u,\sigma\right)  $ such that for any $\varepsilon\in\left(  0,1\right]  $,
\[
\left\|  u\right\|  _{\infty}^{2}\leq\exp\left(  e^{C/\varepsilon}\right)
\text{ ,}%
\]

\noindent which allows us to get to the existence of a fixed point for
$\mathbf{H}\left(  \cdot,1\right)  $. Furthermore, by (\ref{2.13}), the
control is then bounded as follows: for any $\varepsilon\in\left(  0,1\right]
$,
\[
\left\|  h\right\|  _{L^{2}\left(  \omega\times\left(  0,T\right)  \right)
}\leq\exp\left(  e^{C/\varepsilon}\right)  \text{ .}%
\]

\noindent This completes the proof.

\bigskip

\bigskip

Remark .- Notice that the measure of the cost of the control of the semilinear
heat equation (\ref{1.2}) can be improved and become of order
$e^{C/\varepsilon^{2}}$ by adding the following more restrictive hypothesis
$f\left(  0\right)  =0$ and
\[
\underset{\left|  s\right|  \rightarrow\infty}{\lim}\frac{f\left(  s\right)
}{\left|  s\right|  \sqrt{\left|  \ln\ln\left(  1+\left|  s\right|  \right)
\right|  }}=0\text{ .}%
\]

\noindent Indeed, the minimization of the second member of (\ref{2.12}) with
respect to the quantity $\left\|  g\left(  u\right)  \right\|  _{\infty}$
suggests us our choice (\ref{2.9}) of the time of controllability $T^{\prime}%
$. But the minimization of the second member of (\ref{2.12}) when $f\left(
0\right)  =0$ with respect to $\varepsilon\in\left(  0,1\right]  $, suggests
to take $T^{\prime}=\varepsilon T$ in order to get an estimate of the cost of
order $e^{C/\varepsilon^{2}}$.

\bigskip

\bigskip
\end{document}